\documentclass[12pt]{article}

\usepackage{amsmath,amsfonts}
\usepackage[english]{babel}
\usepackage[latin1]{inputenc}
\usepackage{hyperref}

\usepackage{geometry} 	
\numberwithin{equation}{section}
\newtheorem{theorem}{\textbf{Theorem}}
\numberwithin{theorem}{section}

\numberwithin{corollary}{section}

\newtheorem{remark}[theorem]{\textbf{Remark}}

\def\QED{{\boldmath$\rule{0.5em}{0.5em}$}}                                
\def\markatright#1{\leavevmode\unskip\nobreak\quad\hspace*{\fill}{#1}}    
\def\qed{\markatright{\QED}}                                              

\usepackage{authblk}

\title{\bf A two-sided Faulhaber-like formula involving Bernoulli polynomials}
\author[1,4]{J. Fernando Barbero G.}
\author[2,4,5]{Juan Margalef-Bentabol}
\author[3,4]{Eduardo J.S. Villase\~nor\,}
\affil[1]{Instituto de Estructura de la Materia, CSIC. Serrano 123, 28006 Madrid, Spain}
\affil[2]{Laboratory of Geometry and Dynamical Systems, Department of Mathematics, EPSEB, Universitat Polit\`ecnica de Catalunya, BGSMath, Barcelona, Spain}
\affil[3]{Universidad Carlos III de Madrid. Avda.\  de la Universidad 30, 28911 Legan\'es, Spain}
\affil[4]{Grupo de Teor\'{\i}as de Campos y F\'{\i}sica Estad\'{\i}stica. Instituto Gregorio Mill\'an (UC3M). Unidad Asociada al Instituto de Estructura de la Materia, CSIC, Madrid, Spain}
\affil[5]{Institute for Gravitation and the Cosmos \& Physics Department, Penn State, University Park, PA 16802, USA}

\date{}                     
\begin{document}
	
	\maketitle
	
	\begin{abstract}
		\noindent We give a new identity involving Bernoulli polynomials and combinatorial numbers. This provides, in particular, a Faulhaber-like formula for sums of the form $\displaystyle 1^m(n-1)^m+2^m(n-2)^m+\cdots+(n-1)^m1^m$ for positive integers $m$ and $n$.
	\end{abstract}

	\section{Introduction} \label{sect1}
	Bernoulli numbers $B_k$ are given by the exponential generating function $z/(e^z-1)$,
	\[B_k=k![z^k]\frac{z}{e^z-1}\,,\]
	where $[z^n]f(z)$ is the $n$-th coefficient of the Taylor expansion of $f$ around $z=0$.\medskip
	
	In the course of studying the distribution of the eigenvalues of the so-called \textit{area operator} in loop quantum gravity \cite{area operator} we were led to believe that the following identity held
	
	\begin{equation}\label{equation: identidad w=-2}
	\sum_{k=0}^m\binom{m}{k}\frac{B_{2m-k+1}}{2m-k+1}2^{k}=\frac{(-1)^m}{2}\left(1-2^{2m+1}\frac{\Gamma(1+m)^2}{\Gamma(2m+2)}\right)\,,
	\end{equation}
	for $m\in\mathbb{N}\cup\{0\}$. The purpose of this short note is to prove this formula by proving a generalization of it. Particular cases of this general formula involve what we called a two-sided Faulhaber-like formula. A Faulhaber formula (also called Bernoulli's formula as Jacob Bernoulli was the first to write it) is given by
	\[\sum_{k=1}^n k^p=\frac{1}{p+1}\sum_{j=0}^p\binom{p+1}{j}B_jn^{p+1-j}.\]
	Notice that in
	\[\sum_{k=1}^{n-1} k^p=1^p+2^p+\cdots+(n-2)^p+(n-1)^p\]
	there is an increasing sequence of addends given by powers of the integers. A particular and interesting case of the aforementioned generalized formula will involve instead a ``two-sided'' version of it:
	\[\sum_{k=1}^{n-1} k^p(n-k)^p=1^p(n-1)^p+2^p(n-2)^p+\cdots+(n-2)^p2^p+(n-1)^p1^p.\]
	
	Likewise, the Bernoulli numbers are generalized by considering the Bernoulli polynomials:
	\[B_k(x)=k![z^k]\frac{z e^{xz}}{e^z-1}\,.\]
	
	\section{Main theorem} \label{sect2}
	
	The main result of the paper is the following
	
	\begin{theorem}\label{th1}
		Given $N\in\mathbb{Z}$, $m\in\mathbb{N}$  and $w\in\mathbb{C}$, we have
		\begin{align}
		&\sum_{k=0}^m\binom{m}{k}\frac{B_{m+k+1}\left(\frac{N-w}{2}\right)}{m+k+1}w^{m-k}\label{equation: identidad general}\\
		&=\frac{(-1)^{m+1}}{2^{2m+1}}\!\left[\!\frac{(2w)^{2m+1}}{2(2m+1)\binom{2m}{m}}\!-\mathrm{sign}(N-1)\!\!\sum_{k=1}^{|N-1|}\!\!\Big(\!w^2\!-\!(|N-1|\!-2k+1)^2\!\Big)^m\!\right].\nonumber
		\end{align}	
	\end{theorem}
	
	Before proceeding with the proof let us discuss some consequences of this formula

	\begin{remark}
		It is possible to get a number of Faulhaber-like formulas from \eqref{equation: identidad general}. The simplest one can be obtained by taking both $w$ and $N$ to be equal to a natural number $n\geq2$.
		\begin{align}\label{Faulhaber}
		\begin{split}
		\sum_{k=1}^{n-1} k^m(n-k)^m&=\frac{n^{2m+1}}{(2m+1)\binom{2m}{m}}+2(-1)^m\sum_{k=0}^m\binom{m}{k}\frac{B_{m+k+1}}{m+k+1}n^{m-k}=\\
		&=\frac{n^{2m+1}}{(2m+1)\binom{2m}{m}}-2(-1)^m\sum_{k=0}^m\binom{m}{k}\zeta(-m-k)n^{m-k}\,.
		\end{split}
		\end{align}
		where we have used the well known relation between the zeta Riemann function and the Bernoulli numbers
		\[\zeta(1-N)=-\frac{B_{N}}{N}.\]
		Equation \eqref{Faulhaber} appears often in the literature obtained through different methods (see for instance \cite[page 10]{kolosov}).
		
	\end{remark}
	
	\begin{remark}
		For $N=1$, Equation \eqref{equation: identidad general} gives the beautiful expression (equivalent to equation (1.17) of \cite{Sun}) valid for any $w\in\mathbb{C}$,
		\begin{equation}\label{Formula_w}
		\sum_{k=0}^m\binom{m}{k}\frac{B_{m+k+1}\left(\frac{1-w}{2}\right)}{m+k+1}w^{m-k}=\frac{(-1)^{m+1}}{2}\frac{w^{2m+1}}{(2m+1)\binom{2m}{m}}\,.
		\end{equation}
	\end{remark}
	
	\begin{remark} Sums involving
		\[\frac{B_{\beta m+k+1}}{\beta m+k+1}=-\zeta(-k-\beta m)\]
		with integer $\beta\geq2$ can also be studied although a more complicated approach is needed involving complex analysis and combinatorial identities. Nonetheless, the results are not as neat as \eqref{equation: identidad general} and each case has to be studied separately.
	\end{remark}
	
	\begin{remark} It is also possible to generalize \eqref{equation: identidad general} for fractional values of $N$ but, again, no systematic approach has been found. One such expression is when $w=N=1/2$
		\[
		(-1)^{m+1}\sum_{k=0}^m\binom{m}{k}\frac{B_{m+k+1}}{m+k+1}2^k=\frac{1}{2^{m+2}(2m+1)\binom{2m}{m}}+\frac{1}{2^{3m+2}}\sum_{k=0}^m(-1)^k\binom{m}{k}E_{2k}
		\]
		where the $E_n$ are the Euler numbers \cite[entry A122045]{OEIS}.
	\end{remark}
	
	\mbox{}
	
	\noindent\textbf{Proof of Theorem \ref{th1}}
	
	\noindent The result is a consequence, on one hand, of the following easy-to-prove formula for the Bernoulli polynomials
	\begin{align}\begin{split}
	&B_n(x+r)-B_n(x)=n\,\mathrm{sign}(r)\left(\sum_{k=1}^{|r|-1}(x+k\,\mathrm{sign}(r)-1)^{n-1}\right.\label{id_Bernoulli}\\
	&\left.\hspace*{1cm}\phantom{\sum_{k=1}^{|r|-1}}+\frac{1+\mathrm{sign}(r)}{2}(x+r-1)^{n-1}+\frac{1-\mathrm{sign}(r)}{2}(x-1)^{n-1}\right)\,,
	\end{split}
	\end{align}
	valid for $r\in\mathbb{Z}$ and $x\in \mathbb{C}$, which is a direct consequence of \[B_n(x+1)-B_n(x)=n x^{n-1}\,,\] and, on the other hand, of the remarkable identity obtained by Sun (equation (1.14) of \cite{Sun})
	\begin{align}\label{equation_Sun}
	\begin{split}&(-1)^k \sum_{j=0}^k\binom{k}{j}x^{k-j}\frac{B_{\ell+j+1}(y)}{\ell+j+1}+(-1)^\ell \sum_{j=0}^\ell\binom{\ell}{j}x^{\ell-j}\frac{B_{k+j+1}(z)}{k+j+1} \\
	&=\frac{(-x)^{k+\ell+1}}{(k+\ell+1)\binom{k+\ell}{k}}
	\end{split}
	\end{align}
	where $k\,,\ell\in\mathbb{N}$ and $x+y+z=1$. \medskip
	
	\noindent Taking now $x=w$, $y=(N-w)/2$, $z=1-(N+w)/2$ and $k=\ell=m\in\mathbb{N}$ in \eqref{equation_Sun} we obtain
	\begin{align*}
	& (-1)^m\sum_{j=0}^m\binom{m}{j}w^{m-j}\frac{B_{m+k+1}\left(\frac{N-w}{2}\right)}{m+j+1} \\
	& =\frac{(-w)^{2m+1}}{(2m+1)\binom{2m}{m}}+ (-1)^{m+1}\sum_{j=0}^m\binom{m}{j}w^{m-j}\frac{B_{m+j+1}\left(1-\frac{N+w}{2}\right)}{m+j+1}\,.
	\end{align*}
	Using now equation \eqref{id_Bernoulli} to rewrite the last term in terms of $B_{m+j+1}\left(\frac{N-w}{2}\right)$, we finally obtain \eqref{equation: identidad general}.\hfill \qed
	
	\bigskip
	
	\section*{Acknowledgments}
	This work has been supported by the Spanish Ministerio de Ciencia Innovaci\'on y Universidades-Agencia Estatal de Investigaci\'on/FIS2017-84440-C2-2-P grant.  Juan Margalef-Bentabol is supported by 2017SGR932 AGAUR/Generalitat de Catalunya, MTM2015-69135-P/FEDER, MTM2015-65715-P, and the ERC Starting Grant with number 335079. He is also supported in part by the Eberly Research Funds of Penn State, by the NSF grant PHY-1806356, and by the Urania Stott fund of Pittsburgh foundation UN2017-92945.

\end{document}